\documentclass[11pt]{article}
\usepackage{bbm}
\usepackage{mathrsfs}
\usepackage{amscd}
\usepackage{amsmath,amsfonts,amssymb,amscd}
\usepackage{indentfirst,graphics,epsfig,psfrag}
\input{epsf}
\usepackage{ifpdf}
\usepackage{enumerate}
\usepackage{appendix}
\usepackage{enumerate}
\usepackage{amsmath}
\usepackage{lineno}
\usepackage{color}
\makeatletter \@addtoreset{figure}{section} \makeatother
\makeatletter
\long\def\@makecaption#1#2{%
   \vskip 10\p@
   \setbox\@tempboxa\hbox{{#1}\ \ #2}%
   \ifdim \wd\@tempboxa >\hsize

       {#1}\ \ #2\par
   \else
       \hbox to\hsize{\hfil\box\@tempboxa\hfil}%
   \fi}
\makeatother

\newtheorem{thm}{Theorem}[section]

\newtheorem{lem}{Lemma}[section]

\newtheorem{obs}[thm]{Observation}
\newtheorem{pro}{Proposition}[section]

\newcommand{\qed}{{\hfill\rule{3pt}{7pt}}}

\def\qed{\hfill \rule{4pt}{7pt}}

\begin{document}
\title{\textbf{On the Equitable Vertex Arboricity of Graphs}
\footnote{Research supported by the National Science Foundation of
China (Nos. 11551001, 61164005, 11101232, 11461054 and 11161037),
the National Basic Research Program of China (No. 2010CB334708) and
the Program for Changjiang Scholars and Innovative Research Team in
Universities (No. IRT1068), the Research Fund for the Chunhui
Program of Ministry of Education of China (No. Z2014022) and the
Nature Science Foundation from Qinghai Province (Nos. 2012-Z-943,
2014-ZJ-721 and 2014-ZJ-907).}}
\author{
\small Yaping Mao$^1$, \ Zhiwei Guo$^{1}$, \ Hong-Jian Lai$^2$\footnote{Corresponding author}, \ Haixing Zhao$^3$\\[0.3cm]
\small $^1$Department of Mathematics, Qinghai Normal\\
\small University, Xining, Qinghai 810008, China\\[0.2cm]
\small $^2$Department of Mathematics, West Virginia\\
\small University, Morgantown, WV 26506, USA\\[0.2cm]
\small $^3$School of Computer, Qinghai Normal\\
\small University, Xining, Qinghai 810008, China\\[0.2cm]
\small E-mails: maoyaping@ymail.com; guozhiweic@sohu.com;\\
\small hjlai2015@hotmail.com; h.x.zhao@163.com}
\date{}
\maketitle
\begin{abstract}
The equitable coloring problem, introduced by Meyer in 1973, has
received considerable attention and research. Recently, Wu, Zhang
and Li introduced the concept of equitable $(t,k)$-tree-coloring,
which can be regarded as a generalization of proper equitable
$t$-coloring. The \emph{strong equitable vertex $k$-arboricity} of
$G$, denoted by ${va_k}^\equiv(G)$, is the smallest integer $t$ such
that $G$ has an equitable $(t', k)$-tree-coloring for every $t'\geq
t$. The exact value of strong equitable vertex $k$-arboricity of
complete equipartition bipartite graph $K_{n,n}$ was studied by Wu,
Zhang and Li. In this paper, we first get a sharp upper bound of
strong equitable vertex arboricity of complete bipartite graph
$K_{n,n+\ell} \ (1\leq \ell\leq n)$, that is,
${va_2}^\equiv(K_{n,n+\ell})\leq
2\left\lfloor{\frac{n+\ell+1}{3}}\right\rfloor$. Next, we obtain a
sufficient and necessary condition on an equitable $(q,
\infty)$-tree coloring of a complete equipartition tripartite graph,
and study the strong equitable vertex arboricity of forests. For a
simple graph $G$ of order $n$, we show that $1\leq
{va_k}^\equiv(G)\leq \lceil n/2 \rceil$. Furthermore, graphs with
${va_k}^\equiv(G)=1,\lceil\frac{n}{2}\rceil,\lceil\frac{n}{2}\rceil-1$
are characterized, respectively. In the end, we obtain the
Nordhaus-Gaddum type results of strong equitable vertex
$k$-arboricity for general $k$.
\\[2mm]
{\bf Keywords:} Equitable coloring, vertex $k$-arboricity,
$k$-tree-coloring, complete multipartite graph, Nordhaus-Gaddum type result.
\\[2mm]
{\bf AMS subject classification 2010:} 05C05; 05C12; 05C35.
\end{abstract}

\section{Introduction}

All graphs considered in this paper are finite and simple. For a
real number $x$, $\lceil x\rceil$ is the least integer not less than
$x$ and $\lfloor x \rfloor$ is the largest integer not larger than
$x$. We use $V(G)$, $E(G)$, $\delta(G)$ and $\Delta(G)$ to denote
the vertex set, edge set, minimum degree and maximum degree of $G$,
respectively. For a vertex $v \in V(G)$, let $ N_G(v)$ denote the
set of neighbors of $v$ in $G$ and $d_G(v)=|N_G(v)|$ denote the
degree of $v$. We often use $d(v)$ for $d_G(v)$ and refer the reader
to \cite{Bollo} for undefined terms and notation.

A mapping $f: V(G) \rightarrow\{1,2,\ldots,t\}$ is a
\emph{$t$-coloring} of a graph $G$. A $t$-coloring of $G$ is proper
if any two adjacent vertices have different colors. For $1 \leq
i\leq t$, let $V_i = \{v\,|\,f(v) = i\}$. A $t$-coloring of a graph
$G$ is said to be \emph{equitable} if $||V_i|-|V_j||\leq 1$ for all
$i$ and $j$, that is, every color class has size
$\lfloor|V(G)|/t\rfloor$ or $\lceil|V(G)|/t\rceil$. A graph $G$ is
said to be \emph{properly equitably $t$-colorable} if $G$ has a
proper equitable $t$-coloring. The smallest number $t$ for which $G$
is properly equitably $t$-colorable is called \emph{the equitable
chromatic number} of $G$, denoted by $\chi^=(G)$.

The equitable coloring problem, introduced by Meyer
\cite{Meyer}, is motivated by a practical application to municipal
garbage collection \cite{Tucker}. In this context, the vertices of
the graph represent garbage collection routes. A pair of vertices
share an edge if the corresponding routes should not be run on the
same day. It is desirable that the number of routes ran on each day
be approximately the same. Therefore, the problem of assigning one
of the six weekly working days to each route reduces to finding a
proper equitable $6$-coloring. For more applications such as
scheduling, constructing timetables and load balance in parallel
memory systems, we refer to
\cite{Baker,Blazewicz,Das,Irani,Kitagawa,Smith}.

Note that a properly equitably $t$-colorable graph may admit no
proper equitable $t'$-colorings for some $t'> t$. For example, the
complete bipartite graph $H := K_{2m+1,2m+1}$ has no proper
equitable $(2m + 1)$-colorings, although it satisfies $\chi^=(H)=2$.
This fact motivates another interesting parameter for proper
equitable coloring. The equitable chromatic threshold of $G$,
denoted by $\chi^\equiv(G)$, is the smallest integer $t$ such that
$G$ has a proper equitable $t'$-coloring for all $t'\geq t$. This
notion was first introduced by Fan et al. in \cite{Fan}.

For a graph $G$, a \emph{$d$-relaxed $k$-coloring}, also known as a
$d$-defective coloring, of $G$ is a function $f$ from $V(G)$ to
$\{1,2,\ldots,k\}$ such that each subgraph $G[V_i]$ is a graph of
maximum degree at most $d$. A \emph{$d$-relaxed equitable
$k$-coloring} is a $d$-relaxed $k$-coloring that is equitable. In
\cite{Fan}, Fan, Kierstead, Liu, Molla, Wu and Zhang first
considered relaxed equitable coloring of graphs. They proved that
every graph has a proper equitable $\Delta(G)$-coloring such that
each color class induces a forest with maximum degree at most one.
On the basis of this research, Wu, Zhang and Li \cite{Wu} introduced
the notion of equitable $(t, k)$-tree-coloring, which can be viewed
as a generalization of proper equitable $t$-coloring.

A \emph{$(t,k)$-tree-coloring} is a $t$-coloring $f$ of a graph $G$
such that each component of $G[V_i]$ is a tree of maximum degree at
most $k$. A $(t, \infty)$-tree-coloring is called a
$t$-tree-coloring for short. An \emph{equitable $(t,
k)$-tree-coloring} is a $(t, k)$-tree-coloring that is equitable.
The \emph{equitable vertex $k$-arboricity} of a graph $G$, denoted
by ${va_k}^=(G)$, is the smallest integer $t$ such that $G$ has an
equitable $(t, k)$-tree-coloring. The \emph{strong equitable vertex
$k$-arboricity} of $G$, denoted by ${va_k}^\equiv(G)$, is the
smallest integer $t$ such that $G$ has an equitable $(t',
k)$-tree-coloring for every $t'\geq t$. It is clear that
${va_0}^=(G) = \chi^=(G)$ and ${va_0}^{\equiv}(G) =
{\chi}^{\equiv}(G) $ for every graph $G$, and ${va_k}^=(G)$ and
${va_k}^\equiv(G)$ may vary a lot.

The following results are immediate.
\begin{obs}\label{obs1-1}
If $H$ is spanning subgraph of $G$, then ${va_k}^\equiv(H)\leq
{va_k}^\equiv(G)$.
\end{obs}
\begin{obs}\label{obs1-2}
Let $G$ be a graph of order $n$. Then $G$ has an equitable $(q,
k)$-tree-coloring for every $\lceil \frac{n}{2}\rceil \leq q\leq n$.
\end{obs}

The exact value of strong equitable vertex $1$-arboricity of
complete equipartition bipartite graph $K_{n,n}$ where $n\equiv 2
({\rm mod} \ 3)$ was studied by Wu, Zhang and Li \cite{Wu}. In
Section $2$, we first get a sharp upper bound of strong equitable
vertex arboricity of complete bipartite graph $K_{n,n+\ell} \ (1\leq
\ell\leq n)$, that is, ${va_2}^\equiv(K_{n,n+\ell})\leq
2\left\lfloor{\frac{n+\ell+1}{3}}\right\rfloor$. We next investigate
the strong equitable vertex $k$-arboricity of forests. Wu, Zhang and
Li \cite{Wu} got a sufficient and necessary condition on an
equitable $(q, \infty)$-tree coloring of a complete equipartition
bipartite graph. At the end of Section 2, we obtain a sufficient and
necessary condition on an equitable $(q, \infty)$-tree coloring of a
complete equipartition tripartite graph.

In Section $3$, we show that $1\leq {va_k}^\equiv(G)\leq \lceil n/2
\rceil$ for a simple graph $G$ of order $n$. Furthermore, graphs
with ${va_k}^\equiv(G) \in \{1,\lceil n/2\rceil,\lceil
n/2\rceil-1\}$ are respectively characterized.

Let $\mathcal {G}(n)$ denote the class of simple graphs of order $n
\ (n\geq 2)$. For $G\in \mathcal {G}(n)$, $\Bar{G}$ denotes the
complement of $G$. Give a  graph parameter $f(G)$ and a positive
integer $n$, the \emph{Nordhaus-Gaddum Problem} is to determine
sharp bounds for $(1)$ $f(G)+f(\Bar{G})$ and $(2)$ $f(G)\cdot
f(\Bar{G})$, as $G$ ranges over the class $\mathcal {G}(n)$, and
characterize the extremal graphs, i.e., graphs that achieve the
bounds. The Nordhaus-Gaddum type relations have received wide
attention; see a survey paper \cite{Aouchiche} by Aouchiche and
Hansen in 2013. In Section $4$, we obtain the Nordhaus-Gaddum type
results of strong equitable vertex $k$-arboricity for general $k$.

\section{Results for some specific graphs}

In this section, we give the strong equitable vertex $k$-arboricity
of a complete bipartite graph, a forest and a complete tripartite
graph.

\subsection{Results for complete bipartite graphs}

Wu, Zhang and Li \cite{Wu} obtain the exact value of strong
equitable vertex $1$-arboricity of complete equipartition bipartite
graph $K_{n,n}$ where $n\equiv 2 ({\rm mod} \ 3)$. In this
subsection, we prove a sharp upper bound for the strong equitable
vertex $2$-arboricity of a complete bipartite graph.
\begin{thm}\label{th2-1}
Let $K_{n,n+\ell} \ (1\leq \ell\leq n)$ be a complete bipartite graph. Then
$$
{va_2}^\equiv(K_{n,n+\ell})\leq 2\left\lfloor
{\frac{n+\ell+1}{3}}\right\rfloor.
$$
Moreover, the bound is sharp.
\end{thm}
\begin{pf}
By definition, to show ${va_2}^\equiv(K_{n,n+\ell})\leq 2\lfloor
{\frac{n+\ell+1}{3}}\rfloor$, it suffices to show that
$K_{n,n+\ell}$ has an equitable $(q,2)$-tree-coloring for every $q$
such that $q\geq 2\lfloor {\frac{n+\ell+1}{3}}\rfloor+1$. Note that
$$
3q-2n\geq 6\left\lfloor {\frac{n+\ell+1}{3}}\right\rfloor+3-2n\geq 6
\left(\frac{n+\ell-1}{3}\right)+3-2n=2\ell+1.
$$
Then $q\geq \frac{2n+2\ell+1}{3}$. Furthermore, $\frac{2n+2\ell+1}{3}\leq q\leq 2n+\ell$.

Let $X,Y$ be the partite sets of $K_{n,n+\ell}$ with $|X| = n$ and
$|Y| = n + \ell$. Let $e=xy$ be an edge of $K_{n,n+\ell}$ with $x\in
X$ and $y\in Y$. Color $x$ and $y$ with $1$, and divide each of
$X\setminus \{x\}$ and $Y\setminus \{y\}$ into $\frac{q-1}{2}$
classes equitably and color the vertices of each class with a color
in $\{2,3,\cdots,q\}$. Note that $|X\setminus \{x\}|=n-1$ and
$|Y\setminus \{y\}|=n+\ell-1$.

If $\frac{2n+2\ell+1}{3}\leq q\leq n$, then
$$
2\leq \frac{n-1}{\frac{q-1}{2}}\leq \frac{n+\ell-1}{\frac{q-1}{2}}\leq 3,
$$
and hence each color class contains two or three vertices. Therefore, the resulting coloring is an equitable $(q,2)$-tree-coloring of $K_{n,n+\ell}$.

If $n+\ell\leq q\leq 2n-1$, then
$$
1\leq \frac{n-1}{\frac{q-1}{2}}\leq \frac{n+\ell-1}{\frac{q-1}{2}}\leq 2,
$$
and hence each color class contains one or two vertices. Thus, the resulting coloring is an equitable $(q,2)$-tree-coloring of $K_{n,n+\ell}$, as desired.

If $n\leq q\leq n+\frac{\ell}{2}$, then
$$
2=\frac{2n+\ell}{n+\frac{\ell}{2}}\leq \frac{2n+\ell}{q}\leq \frac{2n+\ell}{n}=2+\frac{\ell}{n}\leq 3.
$$
Hence, each color class contains one or two vertices. Thus, the
resulting coloring is an equitable $(q,2)$-tree-coloring of
$K_{n,n+\ell}$, as desired.

If $n+\frac{\ell}{2}+1\leq q\leq n+\ell$, then
$$
1\leq 1+\frac{n}{n+\ell}=\frac{2n+\ell}{n+\ell}\leq \frac{2n+\ell}{q}\leq \frac{2n+\ell}{n+\frac{\ell}{2}+1}=2-\frac{4}{2n+\ell-2}\leq 2.
$$
Hence, each color class contains one or two vertices. Thus, the
resulting coloring is an equitable $(q,2)$-tree-coloring of
$K_{n,n+\ell}$, as desired.

If $2n\leq q\leq 2n+\ell$, then
$$
1=\frac{2n+\ell}{2n+\ell}\leq \frac{2n+\ell}{q}\leq \frac{2n+\ell}{2n}=1+\frac{\ell}{2n}\leq 2.
$$
Hence, each color class contains one or two vertices. Thus, the resulting coloring is an equitable $(q,2)$-tree-coloring of $K_{n,n+\ell}$, as desired.\qed
\end{pf}

\begin{lem}\label{lem2-1}
If the cardinality of color classes is at least $4$ in any equitable
$(t,2)$-tree-coloring of $K_{n,n+\ell}$, then all vertices of each
color class must belong to the same partite set.
\end{lem}
\begin{pf}
Assume, to the contrary, that all the vertices in some color class
belong to different partite sets of $K_{n,n+\ell}$. Then as this
color class has at least 4 vertices, the subgraph induced by all
vertices of this color class either has a cycle, contrary to that
this subgraph should be a tree; or has maximum degree at least $3$.
Contrary to the definition of equitable $(t,2)$- tree-colorings. So,
the conclusion holds.\qed
\end{pf}
\begin{pro}\label{pro2-1}
Let $n=3t \ (t\leq 2)$. Then ${va_2}^\equiv(K_{n,n+1})=
2\left\lfloor {\frac{n+2}{3}}\right\rfloor.$
\end{pro}
\begin{pf}
Suppose $n=3t \ (t\leq 2)$ and $\ell=1$. By Theorem \ref{th2-1}, we
have ${va_2}^\equiv(K_{n,n+1})\leq 2\left\lfloor
{\frac{n+2}{3}}\right\rfloor$. Next, in order to prove that
${va_2}^\equiv(K_{n,n+1})\geq 2\left\lfloor
{\frac{n+2}{3}}\right\rfloor$, we need to prove that $K_{n,n+1}$ has
no equitable $(2t-1, 2)$-tree-coloring. Assume, to the contrary,
that $K_{n,n+1}$ has an equitable $(2t-1, 2)$-tree-coloring. Then
cardinality of each color class is at least $4$ because
$\lfloor{\frac{2n+1}{2t-1}}\rfloor = \lfloor
{\frac{6t+1}{2t-1}}\rfloor=4$ for $t\leq2$. By Lemma \ref{lem2-1},
we have that all the vertices of every color class belong to some
partite set of $K_{n,n+1}$. However, for $t\leq2$, $n$ and $n+1$
both are divisible by $4$. a contradiction. So,
${va_2}^\equiv(K_{n,n+1})= 2\left\lfloor
{\frac{n+2}{3}}\right\rfloor$. \qed
\end{pf}

From Proposition \ref{pro2-1}, we can see that this bound of Theorem
\ref{th2-1} is sharp.

\subsection{Results for forests}

We first give the exact value of strong equitable vertex
$k$-arboricity of a wheel, which also be used later.
\begin{lem}\label{lem2-2}
Let $W_{n}$ be a wheel of order $n$. Then
$$
{va_k}^\equiv(W_{n})\leq \left\lceil\frac{n}{k}\right\rceil.
$$
Moreover, the bound is sharp.
\end{lem}
\begin{pf}
From Observation \ref{obs1-2}, we only need to prove that $W_{n}$
has an equitable $(q,k)$-tree-coloring for each $q$ with
$\lceil\frac{n}{k}\rceil\leq q\leq \lceil n/2\rceil$. Set
$V(W_{n})=\{v_1,v_2,\ldots,v_n\}$. Without loss of generality, let
$v_1$ be the center of $W_{n}$. For $1\leq i\leq q$, we set
$V_i=\{v_{jq+i}\in V(W_{n})\,|\,0\leq j\leq k\}$. We give a
vertex-coloring $c$ of $V(W_{n})$ with $q$ colors such that all the
vertices of $V_i$ receive the color $i \ (1\leq i\leq q)$. Clearly,
the maximum degree of the induced subgraph of each $V_i$ is at most
$k$. So ${va_k}^\equiv(W_{n})\leq \lceil n/k\rceil$, as desired.
\end{pf}

\begin{pro}\label{pro}
Let $K_{1,n-1}$ be a star of order $n$. Then
$$
1+\left\lceil {\frac{n-k-1}{k+2}}\right\rceil\leq
{va_k}^\equiv(K_{1,n-1})\leq \left\lceil\frac{n}{k}\right\rceil.
$$
Moreover, the bound is sharp.
\end{pro}
\begin{pf}
From Lemma \ref{lem2-2} and Observation \ref{obs1-1},
${va_k}^\equiv(K_{1,n-1})\leq \left\lceil \frac{n}{k} \right\rceil$.
We now show that ${va_k}^\equiv(K_{1,n-1})\geq 1+\lceil
{\frac{n-k-1}{k+2}}\rceil$. In order to show this, we need to prove
that $K_{1,n-1}$ has no equitable $(\lceil
{\frac{n-k-1}{k+2}}\rceil,k)$-tree-coloring. Assume, to the
contrary, that $K_{1,n-1}$ has an equitable $(\lceil
{\frac{n-k-1}{k+2}}\rceil,k)$-tree-coloring. Since
$$
\frac{n}{\left\lceil\frac{n-k-1}{k+2}\right\rceil}\geq
\frac{n}{\frac{n}{k+2}}=k+2,
$$
it follows that each color class contains at least $k+2$ vertices.
Then the color class containing the center of $K_{1,n-1}$ also
contains at least $k+2$ vertices, and hence the maximum degree of
the subgraph induced by all the vertices of this color class is at
least $k+1$, a contradiction. So ${va_k}^\equiv(K_{1,n-1})\geq
1+\lceil {\frac{n-k-1}{k+2}}\rceil$.
\end{pf}

\begin{pro}\label{pro2-3}
Let $G$ be a forest with maximum degree $\Delta$. If $k\geq \Delta$,
then ${va_k}^\equiv(G)=1$.
\end{pro}
\begin{pf}
Suppose that $G$ is a forest with the maximum degree $\Delta$ and
$k\geq \Delta$. The graph induced by all vertices is a forest of
maximum degree at most $\Delta$. From the definition of the strong
equitable vertex $k$-arboricity of $G$, we have
${va_\Delta}^\equiv(G)=1$.\qed
\end{pf}

\subsection{Results for complete tripartite graphs}

In this section, we will prove a sufficient and necessary condition
for a partial $q$-coloring of $K_{n,n,n}$ to be an equitable $(q,
\infty)$-tree coloring. Let $K_{n,n,n}$ be a complete tripartite
graph with three partite sets $X$, $Y$ and $Z$. For a partial
$q$-coloring $c$ of $K_{n,n,n}$ (not need to be proper), we let $V_i
\ (1\leq i\leq q$) be its color classes and
$a=\lfloor\frac{3n}{q}\rfloor$. Set
\begin{itemize}
\item $X_1=\{V_i\ |V_i\cap X|=a+1,|V_i\cap Y|=0,|V_i\cap Z|=0\}$;

\item $X_1'=\{V_i|\ |V_i\cap X|=a,|V_i\cap Y|=1,|V_i\cap Z|=0\}$;

\item $X_1''=\{V_i|\ |V_i\cap X|=a,|V_i\cap Y|=0,|V_i\cap Z|=1\}$;

\item $X_2=\{V_i\ |V_i\cap X|=a,|V_i\cap Y|=0,|V_i\cap Z|=0\}$;

\item $X_2'=\{V_i|\ |V_i\cap X|=a-1,|V_i\cap Y|=1,|V_i\cap Z|=0\}$;

\item $X_2''=\{V_i|\ |V_i\cap X|=a-1,|V_i\cap Y|=0,|V_i\cap Z|=1\}$;

\item $Y_1=\{V_i\ |V_i\cap Y|=a+1,|V_i\cap X|=0,|V_i\cap Z|=0\}$;

\item $Y_1'=\{V_i|\ |V_i\cap Y|=a,|V_i\cap X|=1,|V_i\cap Z|=0\}$;

\item $Y_1''=\{V_i|\ |V_i\cap Y|=a,|V_i\cap X|=0,|V_i\cap Z|=1\}$;

\item $Y_2=\{V_i\ |V_i\cap Y|=a,|V_i\cap X|=0,|V_i\cap Z|=0\}$;

\item $Y_2'=\{V_i|\ |V_i\cap Y|=a-1,|V_i\cap X|=1,|V_i\cap Z|=0\}$;

\item $Y_2''=\{V_i|\ |V_i\cap Y|=a-1,|V_i\cap X|=0,|V_i\cap Z|=1\}$;

\item $Z_1=\{V_i\ |V_i\cap Z|=a+1,|V_i\cap X|=0,|V_i\cap Y|=0\}$;

\item $Z_1'=\{V_i|\ |V_i\cap Z|=a,|V_i\cap X|=1,|V_i\cap Y|=0\}$;

\item $Z_1''=\{V_i|\ |V_i\cap Z|=a,|V_i\cap X|=0,|V_i\cap Y|=1\}$;

\item $Z_2=\{V_i\ |V_i\cap Z|=a,|V_i\cap X|=0,|V_i\cap Y|=0\}$;

\item $Z_2'=\{V_i|\ |V_i\cap Z|=a-1,|V_i\cap X|=1,|V_i\cap Y|=0\}$;

\item $Z_2''=\{V_i|\ |V_i\cap Z|=a-1,|V_i\cap X|=0,|V_i\cap Y|=1\}$.
\end{itemize}

We now in a position to give our main result.

\begin{thm}\label{th2-2}
If $K_{n,n,n}$ is a complete tripartite graph with three partite
sets $X$, $Y$ and $Z$, where $3n=aq+r$ and $0\leq r\leq q-1$, and
$c$ is a partial $q$-coloring of $K_{n,n,n}$, then $c$ is an
equitable $(q, \infty)$-tree coloring of $K_{n,n,n}$ if and only if
$$
(a+1)|X_1|+a|X_2|+a|X_1'|+a|X_1''|+(a-1)|X_2'|+(a-1)|X_2''|+|Y_1'|+|Y_2'|+|Z_1'|+
|Z_2'|=n \eqno (1)
$$
$$
(a+1)|Y_1|+a|Y_2|+a|Y_1'|+a|Y_1''|+(a-1)|Y_2'|+(a-1)|Y_2''|+|X_1'|+|X_2'|+|Z_1''|+
|Z_2''|=n \eqno (2)
$$
$$
(a+1)|Z_1|+a|Z_2|+a|Z_1'|+a|Z_1''|+(a-1)|Z_2'|+(a-1)|Z_2''|+|X_1''|+|X_2''|+|Y_1''|+
|Y_2''|=n \eqno (3)
$$
\end{thm}
\begin{pf}
Let $V_i \ (1\leq i\leq q$) be the color classes of $c$. Firstly, we
suppose that $c$ is an equitable $(q, \infty)$-tree coloring of
$K_{n,n,n}$. Since $3n=aq+r$ and $0\leq r\leq q-1$, the size of each
color class of $c$ is either $a$ or $a+1$. It is obvious that
$\min\{|V_i\cap X|, |V_i\cap Y|, |V_i\cap Z|\}\leq1$ for every
$1\leq i\leq q$, because otherwise we would find a cycle in some
color class $V_i$, which contradicts to the definition of the
equitable $(t,\infty)$-tree-coloring of $K_{n,n,n,}$. Thus,
$$
\bigcup_{j=1}^2(X_j\cup X_j'\cup X_j''\cup Y_j\cup Y_j'\cup Y_j''\cup Z_j\cup Z_j'\cup Z_j'')=\bigcup_{i=1}^q V_i,
$$
and Equality $(1)\sim(3)$ hold accordingly. On the other hand, if
Equality $(1)\sim(3)$ hold, then $c$ is a $q$-coloring of
$K_{n,n,n}$ and the size of each color class of $c$ is either $a$ or
$a+1$. Furthermore, we also have $\min\{|V_i\cap X|, |V_i\cap Y|,
|V_i\cap Z|\}\leq1$ for every $1\leq i\leq q$. Hence $c$ is an
equitable $(t,\infty)$- tree-coloring of $K_{n,n,n}$. \qed
\end{pf}

\section{Graphs with given strong equitable vertex $k$-arboricity}

In this section, we give the lower and upper bounds for the strong
equitable vertex $k$-arboricity of simple graphs of order $n$.
\begin{pro}\label{pro3-1}
Let $G$ be a simple graph of order $n$. Then
$$
1\leq {va_k}^\equiv(G)\leq \lceil n/2 \rceil.
$$
\end{pro}
\begin{pf}
It is clear that ${va_k}^\equiv(G)\geq 1$. In order to show
${va_k}^\equiv(G)\leq \lceil n/2 \rceil$, we only need to prove that
$G$ has an equitable $(q, k)$-tree-coloring for every $q$ satisfying
 $\lceil n/2 \rceil\leq q\leq n $. If $q=\lceil n/2 \rceil$, then each color class of the resulting equitable tree-coloring of $G$ contains $1$ or $2$ vertices. Suppose $\lceil n/2 \rceil< q\leq n$. We can easily construct an equitable $(q, k)$-tree-coloring of $G$ by coloring the color class of $2$ vertices with two different colors. Hence, we have ${va_k}^\equiv(G)\leq \lceil n/2 \rceil$.\qed
\end{pf}

Graphs with the strong equitable vertex $k$-arboricity can be $1$ and $\lceil\frac{n}{2}\rceil$ are characterized, respectively.

\begin{pro}\label{pro3-2}
Let $G$ be a graph with maximum degree $\Delta$. Then
${va_k}^\equiv(G)=1$ if and only if $k\geq \Delta$ and $G$ is a
forest with the maximum degree $\Delta$.
\end{pro}
\begin{pf}
Suppose ${va_k}^\equiv(G)=1$. By the definition of the strong
equitable vertex $k$-arboricity, graph $G$ induced by all vertices
is a forest such that $\Delta\leq k$. Conversely, assume that $G$ is
a forest with the maximum degree $\Delta$ and $k\geq \Delta$. From
Proposition \ref{pro2-3}, we have ${va_\Delta}^\equiv(G)=1$.\qed
\end{pf}

\begin{thm}\label{th3-3}
Let $G$ be a simple graph of order $n$. Then
${va_k}^\equiv(G)=\lceil\frac{n}{2}\rceil$ if and only if either
$n\geq 3$ is odd and $G=K_n$, or $n$ is even and $G$ satisfies one
of the following conditions.

$(1)$ For $n=2$, $G=K_2$ or $G=2K_1$.

$(2)$ For $n=4$, $\bar{G}$ does not contain $P_4$ as its subgraph.

$(3)$ For $n\geq 6$, every maximum matching $M$ in $\overline{G}$
satisfies $|M|\leq m-1$, where $n=3m+2r$ and $m+r=\lceil n/2
\rceil-1$.
\end{thm}

We give the proof of Theorem \ref{th3-1} by the following lemmas.
\begin{lem}\label{th3-1}
Let $G$ be a simple graph of order $n$, and let $n\geq 3$ be an odd
integer. Then ${va_k}^\equiv(G)=\lceil\frac{n}{2}\rceil$ if and only
if $G=K_n$.
\end{lem}
\begin{pf}
We first consider the case that $n$ is odd. Suppose
${va_k}^\equiv(G)=\lceil\frac{n}{2}\rceil$. We claim that $G=K_n$.
Assume, to the contrary, that $G\neq K_n$. From Observation
\ref{obs1-1}, we only need to show ${va_k}^\equiv(G)\leq
\lceil\frac{n}{2}\rceil-1$ for $G=K_n\setminus e$, where $e\in
E(G)$. Observe that there exists a color class with $3$ vertices,
say $C_1$. Let $n=2\ell+1$. Then $\lceil n/2\rceil=\ell+1$ and
$\frac{n-3}{2}=\frac{2\ell-2}{2}=\ell-1=\lceil n/2\rceil-2$. From
this argument, there are $\lceil n/2\rceil-2$ color classes with $2$
vertices and one color class with $3$ vertices, and hence $G$ has an
equitable $(\lceil n/2\rceil-1,k)$-tree-coloring, which contradicts
to the fact that ${va_k}^\equiv(G)=\lceil\frac{n}{2}\rceil$. So
$G=K_n$.

Conversely, we assume that $G=K_n$. From Proposition \ref{pro3-1},
we have ${va_k}^\equiv(G)\leq \lceil\frac{n}{2}\rceil$. It suffices
to show that ${va_k}^\equiv(G)\geq \lceil\frac{n}{2}\rceil$. Assume,
to the contrary, that $G$ has an equitable $(\lceil
n/2\rceil-1,k)$-tree-coloring. Let $n=2\ell+1 \ (\ell\geq 1)$. Then
$\lceil n/2\rceil-1=\ell$, and
$$
\frac{n}{\lceil n/2\rceil-1}=\frac{2\ell+1}{\ell}=2+\frac{1}{\ell}.
$$
Clearly, we have
$$
2<\frac{n}{\lceil n/2\rceil-1}\leq 3.
$$
Then there exists a color class with $3$ vertices, and the subgraph
induced by the vertices in $C_1$ contains a cycle, a contradiction.
So ${va_k}^\equiv(G)=\lceil\frac{n}{2}\rceil$.\qed
\end{pf}

\begin{lem}\label{lem3-2}
Let $G$ be a simple graph of order $n$, let $n\geq 3$ be an even
integer. Then ${va_k}^\equiv(G)=\lceil n/2\rceil$ if and only if $G$
satisfies one of the following conditions.

$(1)$ For $n=2$, $G=K_2$ or $G=2K_1$.

$(2)$ For $n=4$, $\bar{G}$ does not contain $P_4$ as its subgraph.

$(3)$ For $n\geq 6$, every maximum matching $M$ in $\overline{G}$
satisfies $1\leq |M|\leq m-1$, where $n=3m+2r$ and $m+r=\lceil
n/2\rceil-1$.
\end{lem}
\begin{pf}
Suppose ${va_k}^\equiv(G)=\lceil\frac{n}{2}\rceil$. Set $n=2\ell$.
We claim that $1\leq |M|\leq m-1$ for any maximum matching $M$ in
$\overline{G}$. Assume, to the contrary, that there exists a maximum
matching such that $|M|\geq m$, where $n=3m+2r$ and $m+r=\lceil n/2
\rceil-1$. Set $M=\{v_1w_1,v_2w_2,\ldots,v_rw_r\}$, where $r\geq m$.
Choose $M'=\{v_1w_1,v_2w_2,\ldots,v_mw_m\}$. Set $V(G)\setminus
(\{v_1,v_2,\ldots,v_m\}\cup
\{w_1,w_2,\ldots,w_m\})=\{u_{2m+1},u_{2m+2},\ldots,u_n\}$. Choose
$S_i=\{v_i,w_i,u_{2m+i}\}$, where $1\leq i\leq m$, and choose
$S_{m+j}=\{u_{3m+2j-1},u_{3m+2j}\}$, where $1\leq j\leq
\frac{n-3m}{2}$. Then $S_1,S_2,\ldots,S_{\frac{n-m}{2}}$ are $\lceil
n/2 \rceil-1$ color classes such that the subgraph induced by each
$S_i$ is a forest of order $2$ or $3$. So $G$ has an equitable
$(\lceil n/2\rceil-1,k)$-tree-coloring, which contradicts to the
fact that ${va_k}^\equiv(G)=\lceil\frac{n}{2}\rceil$.

Conversely, we suppose that if any maximum matching $M$ in
$\overline{G}$, then $1\leq |M|\leq m-1$, where $n=3m+2r$ and
$m+r=\lceil n/2 \rceil-1$. From Proposition \ref{pro3-1}, we have
${va_k}^\equiv(G)\leq \lceil\frac{n}{2}\rceil$. It suffices to show
that ${va_k}^\equiv(G)\geq \lceil\frac{n}{2}\rceil$. Assume, to the
contrary, that $G$ has an equitable $(\lceil
n/2\rceil-1,k)$-tree-coloring. Set $n=2\ell \ (\ell\geq 3)$. Then
$\lceil\frac{n}{2}\rceil-1=\ell-1$, and hence
$$
\frac{n}{\lceil\frac{n}{2}\rceil-1}=\frac{2\ell}{\ell-1}=2+\frac{2}{\ell-1}.
$$
Since $\ell\geq 3$, it follows that
$\frac{n}{\lceil\frac{n}{2}\rceil-1}\leq 3$. Let $m$ denote the
number of those color classes such that each color class contains
exactly $3$ vertices, and let $r$ denote the number of those color
classes such that each color class contains exactly $2$ vertices.
Then $n=3m+2r$ and $m+r=\lceil n/2 \rceil-1$. Since $|M|\leq m-1$
for any maximum matching $M$ in $\overline{G}$, there exists a color
class with $3$ vertices, say $C_1$, such that $M\cap
E(\overline{G}[C_1])=\emptyset$. We claim that
$E(\overline{G}[C_1])=\emptyset$. Assume, to the contrary, that
$E(\overline{G}[C_1])\neq \emptyset$. Let $e\in
E(\overline{G}[C_1])$. Then $M\cup \{e\}$ is a matching in
$\overline{G}$, which contradicts to the fact that $M$ is maximum
matching in $\overline{G}$. Therefore,
$E(\overline{G}[C_1])=\emptyset$, and hence $G[C_1]$ is a cycle, a
contradiction. So ${va_k}^\equiv(G)=\lceil\frac{n}{2}\rceil$.\qed
\end{pf}

Graphs with ${va_k}^\equiv(G)=\lceil\frac{n}{2}\rceil-1$ can be also
characterized.
\begin{thm}\label{th3-4}
Let $G$ be a simple graph of order $n \ (n\geq 9, \ n\neq 10)$. Then
${va_k}^\equiv(G)=\lceil\frac{n}{2}\rceil-1$ if and only if $G$
satisfies one of the following conditions. $(1)$ $n\geq 12$ is even,
and every maximum matching $M$ in $\overline{G}$ satisfies $2\leq
|M|\leq m-1$, where $n=3m+2r$ and $m+r=\lceil n/2 \rceil-2$. $(2)$
$n\geq 9$ is odd, and every maximum matching $M$ in $\overline{G}$
satisfies $1\leq |M|\leq m-1$, where $n=3m+2r$ and $m+r=\lceil n/2
\rceil-2$.
\end{thm}
\begin{pf}
We distinguish the following two cases to show our proof.

\textbf{Case 1. $n=2\ell$}

Suppose ${va_k}^\equiv(G)=\lceil\frac{n}{2}\rceil-1$. Then $G$ has
an equitable $(q,k)$-tree-coloring for $q\geq \lceil n/2\rceil-1$.
For $q=\lceil\frac{n}{2}\rceil-1$, since
$$
2<\frac{n}{\lceil\frac{n}{2}\rceil-1}=\frac{2\ell}{\ell-1}\leq 3,
$$
it follows that there is at least one color class of order $3$ in
$G$. We claim that $2\leq |M|\leq m-1$ for any maximum matching $M$
in $\overline{G}$. Assume, to the contrary, that there exists a
maximum matching such that $|M|=1$ or $|M|\geq m$, where $n=3m+2r$
and $m+r=\lceil n/2 \rceil-2=\ell-2$.

Suppose $|M|=1$. Recall that there is at least one color class of
order $3$ in $G$. We furthermore claim that there is only one color
class of order $3$ in $G$. Assume, to the contrary, that there are
two color classes of order $3$ in $G$, say $C_1,C_2$. From the
definition, there is an edge of $C_i \ (i=1,2)$ in $\bar{G}$, say
$e_i$. Clearly, $\{e_1,e_2\}$ is a matching in $\bar{G}$, which
contradicts the fact that $M$ is a maximum matching $\bar{G}$ and
$|M|=1$. So there is only one color class of order $3$ in $G$. Then
each of the remaining classes contains $2$ vertices in $G$. This
contradicts to the fact that $n$ is even.

Suppose that $|M|\geq m$, where $n=3m+2r$ and $m+r=\lceil n/2
\rceil-2=\ell-2$. Set $M=\{v_1w_1,v_2w_2,\ldots,v_rw_r\}$, where
$r\geq m$. Choose $M'=\{v_1w_1,v_2w_2,\ldots,v_mw_m\}$. Set
$V(G)\setminus (\{v_1,v_2,\ldots,v_m\}\cup
\{w_1,w_2,\ldots,w_m\})=\{u_{2m+1},u_{2m+2},\ldots,u_n\}$. Choose
$S_i=\{v_i,w_i,u_{2m+i}\}$, where $1\leq i\leq m$, and choose
$S_{m+j}=\{u_{3m+2j-1},u_{3m+2j}\}$, where $1\leq j\leq
\frac{n-3m}{2}$. Then $S_1,S_2,\ldots,S_{\frac{n-m}{2}}$ are $\lceil
n/2 \rceil-2$ color classes such that the subgraph induced by each
$S_i$ is a forest of order $2$ or $3$. So $G$ has an equitable
$(\lceil n/2\rceil-2,k)$-tree-coloring, which contradicts to the
fact that ${va_k}^\equiv(G)=\lceil\frac{n}{2}\rceil-1$. So $|M|\leq
m-1$.

Conversely, we suppose that every maximum matching $M$ in
$\overline{G}$ satisfies $2\leq |M|\leq m-1$, where $2\leq |M|\leq
m-1$, where $n=3m+2r$ and $m+r=\lceil n/2 \rceil-2$. From Theorem
\ref{th3-3} and Proposition \ref{pro3-1}, we have
${va_k}^\equiv(G)\leq \lceil\frac{n}{2}\rceil-1$. It suffices to
show that ${va_k}^\equiv(G)\geq \lceil\frac{n}{2}\rceil-1$. Assume,
to the contrary, that $G$ has an equitable $(\lceil
n/2\rceil-2,k)$-tree-coloring. Set $n=2\ell \ (\ell\geq 6)$. Then
$\lceil\frac{n}{2}\rceil-2=\ell-2$, and hence
$$
\frac{n}{\lceil\frac{n}{2}\rceil-2}=\frac{2\ell}{\ell-2}=2+\frac{4}{\ell-2}.
$$
Since $\ell\geq 6$, it follows that
$\frac{n}{\lceil\frac{n}{2}\rceil-2}\leq 3$. Let $m$ denote the
number of those color classes such that each color class contains
exactly $3$ vertices, and let $r$ denote the number of those color
classes such that each color class contains exactly $2$ vertices.
Then $n=3m+2r$ and $m+r=\lceil n/2 \rceil-2$. Since $|M|\leq m-1$
for any maximum matching $M$ in $\overline{G}$, there exists a color
class with $3$ vertices, say $C_1$, such that $M\cap
E(\overline{G}[C_1])=\emptyset$. We claim that
$E(\overline{G}[C_1])=\emptyset$. Assume, to the contrary, that
$E(\overline{G}[C_1])\neq \emptyset$. Let $e\in
E(\overline{G}[C_1])$. Then $M\cup \{e\}$ is a matching in
$\overline{G}$, which contradicts to the fact that $M$ is maximum
matching in $\overline{G}$. Therefore,
$E(\overline{G}[C_1])=\emptyset$, and hence $G[C_1]$ is a cycle, a
contradiction. So ${va_k}^\equiv(G)=\lceil\frac{n}{2}\rceil-1$.

\textbf{Case 2. $n=2\ell+1$}

We claim that $1\leq |M|\leq m-1$ for any maximum matching $M$ in
$\overline{G}$. Assume, to the contrary, that there exists a maximum
matching such that $|M|\geq m$, where $n=3m+2r$ and $m+r=\lceil n/2
\rceil-2=\ell-2$. Set $M=\{v_1w_1,v_2w_2,\ldots,v_rw_r\}$, where
$r\geq m$. Choose $M'=\{v_1w_1,v_2w_2,\ldots,v_mw_m\}$. Set
$V(G)\setminus (\{v_1,v_2,\ldots,v_m\}\cup
\{w_1,w_2,\ldots,w_m\})=\{u_{2m+1},u_{2m+2},\ldots,u_n\}$. Choose
$S_i=\{v_i,w_i,u_{2m+i}\}$, where $1\leq i\leq m$, and choose
$S_{m+j}=\{u_{3m+2j-1},u_{3m+2j}\}$, where $1\leq j\leq
\frac{n-3m}{2}$. Then $S_1,S_2,\ldots,S_{\frac{n-m}{2}}$ are $\lceil
n/2 \rceil-2$ color classes such that the subgraph induced by each
$S_i$ is a forest of order $2$ or $3$. So $G$ has an equitable
$(\lceil n/2\rceil-2,k)$-tree-coloring, which contradicts to the
fact that ${va_k}^\equiv(G)=\lceil\frac{n}{2}\rceil-1$. So $1\leq
|M|\leq m-1$.

Conversely, we suppose that if any maximum matching $M$ in
$\overline{G}$ satisfies $1\leq |M|\leq m-1$, where $n=3m+2r$ and
$m+r=\lceil n/2 \rceil-2$. From Theorem \ref{th3-3} and Proposition
\ref{pro3-1}, we have ${va_k}^\equiv(G)\leq
\lceil\frac{n}{2}\rceil-1$. It suffices to show that
${va_k}^\equiv(G)\geq \lceil\frac{n}{2}\rceil-1$. Assume, to the
contrary, that $G$ has an equitable $(\lceil
n/2\rceil-2,k)$-tree-coloring. Set $n=2\ell+1 \ (\ell\geq 4)$. Then
$\lceil\frac{n}{2}\rceil-2=\ell-1$, and hence
$$
\frac{n}{\lceil\frac{n}{2}\rceil-1}=\frac{2\ell+1}{\ell-1}=2+\frac{3}{\ell-1}.
$$
Since $\ell\geq 4$, it follows that
$\frac{n}{\lceil\frac{n}{2}\rceil-2}\leq 3$. Let $m$ denote the
number of those color classes such that each color class contains
exactly $3$ vertices, and let $r$ denote the number of those color
classes such that each color class contains exactly $2$ vertices.
Then $n=3m+2r$ and $m+r=\lceil n/2 \rceil-2$. Since $|M|\leq m-1$
for any maximum matching $M$ in $\overline{G}$, there exists a color
class with $3$ vertices, say $C_1$, such that $M\cap
E(\overline{G}[C_1])=\emptyset$. We claim that
$E(\overline{G}[C_1])=\emptyset$. Assume, to the contrary, that
$E(\overline{G}[C_1])\neq \emptyset$. Let $e\in
E(\overline{G}[C_1])$. Then $M\cup \{e\}$ is a matching in
$\overline{G}$, which contradicts to the fact that $M$ is maximum
matching in $\overline{G}$. Therefore,
$E(\overline{G}[C_1])=\emptyset$, and hence $G[C_1]$ is a cycle, a
contradiction. So ${va_k}^\equiv(G)=\lceil\frac{n}{2}\rceil-1$.\qed
\end{pf}

\section{Nordhaus-Gaddum-type results}

In this section, we investigate the Nordhaus-Gaddum-type problem on
the strong equitable vertex $k$-arboricity of graphs.
\begin{pro}\label{th4-1}
For any $G\in \mathcal {G}(n)$, if $n \ (n\geq 2)$ is even, then

$(1)$ $2\leq {va_k}^\equiv(G)+{va_k}^\equiv(\Bar{G})\leq
2\lceil\frac{n}{2}\rceil$;

$(2)$ $1\leq{va_k}^\equiv(G)\cdot {va_k}^\equiv(\Bar{G})\leq
\lceil\frac{n}{2}\rceil^2$.

Moreover, the bounds are sharp.
\end{pro}
\begin{pf}
$(1)$ From Proposition \ref{pro3-1}, ${va_k}^\equiv(G)\geq 1$ and
${va_k}^\equiv(\Bar{G})\geq 1$, and hence
${va_k}^\equiv(G)+{va_k}^\equiv(\Bar{G})\geq 2$. From Proposition
\ref{pro3-1}, ${va_k}^\equiv(G)\leq \lceil\frac{n}{2}\rceil$.
Furthermore, we have ${va_k}^\equiv(\Bar{G})\leq
\lceil\frac{n}{2}\rceil$ for $\Bar{G}$. Hence,
${va_k}^\equiv(G)+{va_k}^\equiv(\Bar{G})\leq
2\lceil\frac{n}{2}\rceil$, as desired.

$(2)$ From Proposition \ref{pro3-1}, ${va_k}^\equiv(G)\geq 1$ and
${va_k}^\equiv(\Bar{G})\geq 1$, and hence ${va_k}^\equiv(G)\cdot
{va_k}^\equiv(\Bar{G})\geq 1$. From Proposition \ref{pro3-1},
${va_k}^\equiv(G)\leq \lceil\frac{n}{2}\rceil$. For $\Bar{G}$, we
have ${va_k}^\equiv(\Bar{G})\leq \lceil\frac{n}{2}\rceil$. So,
${va_k}^\equiv(G)\cdot {va_k}^\equiv(\Bar{G})\leq
\lceil\frac{n}{2}\rceil^2$. \qed
\end{pf}

\begin{pro}\label{pro5-2}
For any $G\in \mathcal {G}(n)$, if $n \ (n\geq 5)$ is odd, then

$(1)$ $2\leq {va_k}^\equiv(G)+{va_k}^\equiv(\Bar{G})\leq
2\lceil\frac{n}{2}\rceil-2$;

$(2)$ $1\leq{va_k}^\equiv(G)\cdot {va_k}^\equiv(\Bar{G})\leq
(\lceil\frac{n}{2}\rceil-1)^2$.
\end{pro}
\begin{pf}
$(1)$ From Proposition \ref{pro3-1}, ${va_k}^\equiv(G)\geq 1$ and
${va_k}^\equiv(\Bar{G})\geq 1$, and hence
${va_k}^\equiv(G)+{va_k}^\equiv(\Bar{G})\geq 2$. From Proposition
\ref{pro3-1}, ${va_k}^\equiv(G)+{va_k}^\equiv(\Bar{G})\leq
2\lceil\frac{n}{2}\rceil$. Suppose that
${va_k}^\equiv(G)+{va_k}^\equiv(G)=2\lceil\frac{n}{2}\rceil$. Then
${va_k}^\equiv(G)= \lceil\frac{n}{2}\rceil$ and
${va_k}^\equiv(\bar{G})= \lceil\frac{n}{2}\rceil$. From Theorem
\ref{th3-3}, both $G$ and $\bar{G}$ are all complete graph of order
$n$, a contradiction. Suppose that
${va_k}^\equiv(G)+{va_k}^\equiv(G)=2\lceil\frac{n}{2}\rceil-1$. Then
without loss of generality, let ${va_k}^\equiv(G)=
\lceil\frac{n}{2}\rceil$ and ${va_k}^\equiv(\bar{G})=
\lceil\frac{n}{2}\rceil-1$. From Theorem \ref{th3-3}, $G$ is a
complete graph of order $n$, and hence $\bar{G}=nK_1$. Clearly,
${va_k}^\equiv(\bar{G})=1$. From this together with
${va_k}^\equiv(\bar{G})= \lceil\frac{n}{2}\rceil-1$, we have
$\lceil\frac{n}{2}\rceil-1=1$ and $n=4$ or $n=3$, a contradiction.
So ${va_k}^\equiv(G)+{va_k}^\equiv(\Bar{G})\leq
2\lceil\frac{n}{2}\rceil-2$.

$(2)$ From Proposition \ref{pro3-1}, ${va_k}^\equiv(G)\geq 1$ and
${va_k}^\equiv(\Bar{G})\geq 1$, and hence ${va_k}^\equiv(G)\cdot
{va_k}^\equiv(\Bar{G})\geq 1$. From Proposition \ref{pro3-1},
${va_k}^\equiv(G)\cdot {va_k}^\equiv(\Bar{G})\leq
\lceil\frac{n}{2}\rceil^2$. Suppose that ${va_k}^\equiv(G)\cdot
{va_k}^\equiv(G)=\lceil\frac{n}{2}\rceil^2$. Then ${va_k}^\equiv(G)=
\lceil\frac{n}{2}\rceil$ and ${va_k}^\equiv(\bar{G})=
\lceil\frac{n}{2}\rceil$. From Theorem \ref{th3-3}, both $G$ and
$\bar{G}$ are all complete graph of order $n$, a contradiction.
Suppose that ${va_k}^\equiv(G)\cdot
{va_k}^\equiv(G)=\lceil\frac{n}{2}\rceil(\lceil\frac{n}{2}\rceil-1)$.
Similarly to the proof of $(1)$, we can get a contradiction. So
${va_k}^\equiv(G)\cdot
{va_k}^\equiv(G)\leq(\lceil\frac{n}{2}\rceil-1)^2$. \qed
\end{pf}

To show the sharpness of the upper bounds of the above theorems, we
consider the following example.

\noindent\textbf{Example 1.} If $n$ is even, then we let $G=P_4$.
Then $\bar{G}=P_4$. For $k=1$, ${va_1}^\equiv(G)=2$ and
${va_1}^\equiv(\Bar{G})=2$. Then
${va_1}^\equiv(G)+{va_1}^\equiv(\Bar{G})=4=2\lceil\frac{n}{2}\rceil$
and ${va_1}^\equiv(G)\cdot
{va_1}^\equiv(\Bar{G})=4=\lceil\frac{n}{2}\rceil^2$.

If $n$ is odd, then we let $G=C_5$. Then $\bar{G}=C_5$. For $k\geq
1$, ${va_1}^\equiv(G)=2$ and ${va_1}^\equiv(\Bar{G})=2$. Then
${va_1}^\equiv(G)+{va_1}^\equiv(\Bar{G})=4=2\lceil\frac{n}{2}\rceil-2$
and ${va_1}^\equiv(G)\cdot
{va_1}^\equiv(\Bar{G})=4=(\lceil\frac{n}{2}\rceil-1)^2$.\\

Graphs attaining the lower bounds of above theorems can be
characterized.
\begin{pro}\label{pro5-1}
Let $G\in \mathcal {G}(n)$, and let $G$ and $\Bar{G}$ are both
connected. Then ${va_k}^\equiv(G)+{va_k}^\equiv(\Bar{G})=2$ or
${va_k}^\equiv(G)\cdot {va_k}^\equiv(\Bar{G})=1$ if and only if
$G=P_4$ for $k\geq 2$, or $G=P_3$ for $k\geq 2$, or $G=P_2\cup K_1$
for $k\geq 2$,  or $G=P_2$ for $k\geq 1$, or $G=2K_1$ for $k\geq 1$.
\end{pro}
\begin{pf}
If ${va_k}^\equiv(G)+{va_k}^\equiv(\Bar{G})=2$ or
${va_k}^\equiv(G)\cdot {va_k}^\equiv(\Bar{G})=1$, then
${va_k}^\equiv(G)=1$ and ${va_k}^\equiv(\Bar{G})=1$. From
Proposition \ref{pro3-2}, we have $G$ and $\Bar{G}$ are both a
forest and $k\geq \max\{\Delta(G), \Delta(\Bar{G})\}$. Therefore, we
have $2(n-1)\geq E(G)+E(\Bar{G})=E(G\cup
\Bar{G})=E(K_n)=\frac{n(n-1)}{2}$, and hence $n\leq 4$. So, $G=P_4$
for $k\geq 2$, or $G=P_3$ for $k\geq 2$, or $G=P_2\cup K_1$ for
$k\geq 2$,  or $G=P_2$ for $k\geq 1$, or $G=2K_1$ for $k\geq 1$.

Conversely, suppose $G$ satisfies the condition of this theorem. One
can easily check that ${va_k}^\equiv(G)=1$ and
${va_k}^\equiv(\Bar{G})=1$. Hence,
${va_k}^\equiv(G)+{va_k}^\equiv(\Bar{G})=2$ and
${va_k}^\equiv(G)\cdot {va_k}^\equiv(\Bar{G})=1$, as desired. \qed
\end{pf}

\end{document}